\documentclass[11pt]{article}
\usepackage{amstex}
\usepackage{amssymb}
\usepackage{amscd}
\usepackage[all]{xy}
\font\tenmsam=msam10
\def\square{\hbox{\tenmsam\char'003}}

\def\hfl#1#2{\smash{\mathop{\hbox to 12mm{\rightarrowfill}}
\limits^{\scriptstyle#1}_{\scriptstyle#2}}}

\newtheorem{theo}{Theorem}[section]
\newtheorem{prop}[theo]{Proposition}

\newtheorem{lemm}[theo]{Lemma}
\newtheorem{defi}[theo]{Definition}
\newtheorem{rema}[theo]{Remark}
\newtheorem{ex}[theo]{Example}

\title{\textsc{ Galois reconstruction of finite\\
quantum groups}}
\date{} 
\author{Julien Bichon}

\makeatletter
\renewcommand{\@makefnmark}{}
\makeatother

\begin{document}
\maketitle

\begin{abstract}
Let $\mathcal C$ be a (small) category and let 
$F : \mathcal C \longrightarrow \mathcal M alg_f$ be a functor, where
$\mathcal M alg_f$ is the category of finite-dimensional measured algebras
over a field $k$ (or Frobenius algebras).
We construct a universal Hopf algebra $A_{aut}(F)$ such that $F$ factorizes 
through a functor 
$\overline F : \mathcal C \longrightarrow \mathcal M coalg_f(A_{aut}(F))$,
where $\mathcal M coalg_f(A_{aut}(F))$ is the category 
of finite-dimensional measured
 $A_{aut}(F)$-comodule algebras.
This general reconstruction result allows us 
to recapture a finite-dimensional Hopf algebra
$A$ from the category $\mathcal M coalg_f(A)$ and the forgetful functor
$\omega : \mathcal M coalg_f(A) \longrightarrow \mathcal M alg_f$: we have
$A \cong A_{aut}(\omega)$. Our universal construction is also done in
a $C^*$-algebra framework, and we get compact quantum groups in the sense of
Woronowicz. 
\end{abstract}

\footnote{1991 AMS Classification: 16W30, 18D10, 20B25, 46L89.}

\footnote{Keywords: Quantum groups, Galois theory, Tannaka duality, Groups acting on sets.}

\section{Introduction}

We begin with the following two observations: 

{\bf (1)}
Let $\mathcal C$ be a (small) category and 
let $F : \mathcal C \longrightarrow Set_f$
be a functor to finite sets. Let $Aut(F)$ be 
the automorphism group of the functor $F$.
Then $F$ factorizes through a functor 
$\overline F : \mathcal C \longrightarrow 
Aut(F)-Set_f$ (the category of finite sets on which $Aut(F)$ 
operates) followed by the forgetful functor: 

$$\xymatrix{
F : \mathcal C \ar[rr] \ar@{-->}[rd]_{\overline F} && Set_f \\
& Aut(F)-Set_f \ar[ur]}$$

{\bf (2)}
Let $\mathcal C = G-Set_f$ for a finite group $G$ and let $\omega = F$  
 be the forgetful functor. 
Then the groups $Aut(\omega)$ and $G$ are isomorphic.

\medskip

These results are easy to prove, but are of fundamental interest.
Indeed they provide the basis for Grothendieck's Galois theory \cite{[G]},
ie the axiomatic characterization of the categories of the form $G-Sets$
($G$ finite or profinite group).

\medskip

The aim of this paper is to describe quantum 
analogues of the results {\bf (1)} and {\bf (2)}.
We adopt the classical philosophy for quantum spaces and quantum groups : 
a quantum space is thought as an algebra 
(playing the role of the function algebra) and a 
quantum group is thought as a Hopf algebra.
An action of a quantum group on a 
quantum space is  a right comodule algebra structure
on the underlying algebra.

\medskip

Let us note that  we have analogues of {\bf (1)} and {\bf (2)} in the 
framework of linear representations of
algebraic groups, and we also have an axiomatic 
theory known as tannakian duality.
This theory was developed in N. Saavedra Rivano's thesis \cite{[SR]}
under the direction of Grothendieck.
The quantum analogue (ie for non-commutative Hopf algebras) of tannakian  
duality was proved ten years ago by 
K.H. Ulbrich (\cite{[U]}, see also the texts \cite{[JS], [Sc]}).

We use non-commutative tannakian duality to prove our results.
We have also been inspired by S. Wang's discovery \cite{[Wa]}
of the non-existence of the 
quantum automorphism group of a semisimple algebra. As in 
\cite{[Wa]}, this led us to consider measured
quantum spaces (measured algebras or Frobenius algebras, see section 2)
rather than quantum spaces. This framework is in fact very natural, 
especially when we have a view towards the example of comodule algebras
over a finite-dimensional Hopf algebra, thanks to the existence of the Haar 
measure (see \cite{[LS]}).
As a special case of our construction, we get algebraic analogues of the
compact quantum groups constructed in \cite{[Wa]} and \cite{[Bi]}.  

In section 2 we fix some notations and give some preliminary results.
In section 3 we describe  the construction of  a universal Hopf 
algebra $A_{aut}(F)$ associated to any functor 
$F : \mathcal C \longrightarrow \mathcal M alg_f$
(the category of finite-dimensional measured algebras over a field $k$), 
such that $F$ factorizes through a functor
$\overline{F} : \mathcal C \longrightarrow \mathcal M coalg_f(A_{aut}(F))$,
where $\mathcal M coalg_f(A_{aut}(F))$ is the category of 
finite-dimensional measured  $A_{aut}(F)$-comodule algebras.
Let us emphasize that absolutely no structure is required on the category 
$\mathcal C$. This procedure allows us to reconstruct 
a finite-dimensional Hopf algebra $A$ from the category 
$\mathcal M coalg_f(A)$ and the forgetful functor
$\omega : \mathcal M coalg_f(A) \longrightarrow \mathcal M alg_f$. 
Indeed we prove that the Hopf algebras 
$A$ and  $A_{aut}(\omega)$ are isomorphic (under a slight assumption 
on the base field).
This result is a generalization of {\bf (2)}, and is an application of
the fundamental theorem of Hopf modules \cite{[LS],[Mo]}.
It is also relied to the duality between Hopf algebras and multiplicative
unitaries (\cite{[BS]}).

In section 4 we (briefly) describe analogues of the results of section 
3 in a $C^*$-algebra framework. We have to use Woronowicz' Tannaka-Krein
duality \cite{[W2]}, and we get compact quantum groups (\cite{[W1],[W3]}).

\section{Preliminaries}

The aim of this section is to fix some notations and to provide the 
definitions and elementary lemmas (most of them probably well known)
needed in the paper. We work over a fixed field $k$.

Let  $A = (A,m,u,\Delta,\varepsilon,S)$ be a Hopf algebra.
 The multiplication will be denoted by  $m$,   $u : \mathbb C
\rightarrow A$ is the unit of $A$, while $\Delta$, $\varepsilon$ and $S$
are respectively the  comultiplication, the counit and the antipode of $A$.
  
The category of right $A$-comodules will be denoted $Co_f(A)$. If 
$V$ is a right $A$-comodule, the coaction will be denoted
by $\alpha_V : V \longrightarrow V \otimes A$. 
The Hopf algebra $A$ is always a right $A$-comodule 
(let us say the regular one) via the comultiplication $\Delta$.

Let $V$ be a right $A$-comodule. We denote by $V_0$ the associated 
trivial comodule: the underlying vector space is $V$ and the coaction
is $1_V \otimes u$. In this way the coaction
becomes a comodule map $\alpha_V : V \longrightarrow V_0 \otimes A$.

The first two lemmas are used in the proof of theorem
3.5.
The first one is a well-known trick in Hopf algebra theory
(see lemma 3.1.4 in \cite{[Mo]}).

\begin{lemm}
Let $A$ be a Hopf algebra and let $V$ be a right $A$-comodule.
Then the map $\kappa_V : V \otimes A \longrightarrow V_0 \otimes A$
defined by the composition
$$
V \otimes A \stackrel{\alpha_V \otimes 1_A}{\longrightarrow}
V_0 \otimes A \otimes A
\stackrel{1_{V_0} \otimes m}{\longrightarrow} V_0 \otimes A$$
is an isomorphism of $A$-comodules. \square 
\end{lemm} 

Indeed, the inverse is given by $(1_V \otimes m) \circ (1_V \otimes S
 \otimes 1_A)
\circ (\alpha_V \otimes 1_A)$. When $V= A$ we have a pentagonal operator
as in \cite{[BS]}. 

Let us recall that an $A$-Hopf module is a vector space $V$ which both 
a right $A$-comodule and a left $A$-module such that the action
$\mu : A \otimes V \longrightarrow V$ is an $A$-comodule map.

\begin{lemm}
Let $H$ and $A$ be Hopf algebras and let 
$\pi : H \longrightarrow A$ be a Hopf algebra morphism.
Let us assume:

\noindent
i) There is a $H$-comodule algebra structure
$\alpha : A \longrightarrow A \otimes H$ on $A$ such that
$(1_A \otimes \pi) \circ \alpha = \Delta_A$.

\noindent
ii) The map $\kappa = \kappa_A = (1_A \otimes m_A) \circ 
(\Delta_A \otimes 1_A) : A\otimes A \longrightarrow A_0 \otimes A$
(where $A_0$ denotes the trivial $H$-comodule structure on $A$) is a
$H$-comodule map.

Then $A$, endowed with its $H$-comodule structure and its left
$H$-module structure induced by $\pi$, is a $H$-Hopf module.
If furthermore the antipode of $H$ is bijective and $A$ is finite-dimensional
then $\pi$ is an isomorphism. 
\end{lemm}

\noindent
{\bf Proof}. Let $\mu = m_A \circ (\pi \otimes 1_A) : H \otimes A \longrightarrow A$. We must show that $\mu$ is a map of $H$-comodules, that is:
$$\alpha \circ \mu = (\mu \otimes 1_H)
\circ (1_H \otimes 1_A \otimes m_H) \circ
(1_H \otimes C_{H,A} \otimes 1_H) \circ (\Delta_H \otimes \alpha) \ (\star)
$$
where $C_{H,A}$ denotes the symmetry. By assumption we have:
\begin{align*}
(1_A & \otimes\alpha)  \circ (1_A \otimes m_A) 
\circ (\Delta_A \otimes 1_A) = \\
& (1_A \otimes m_A \otimes 1_H) \circ 
(\Delta_A \otimes 1_A \otimes 1_H) \circ
(1_A \otimes 1_A \otimes m_H) \circ
(1_A \otimes C_{H,A} \otimes 1_H) \circ (\alpha \otimes \alpha)
\end{align*}
and since $\Delta_A = (1_A \otimes \pi) \circ \alpha$, we have
$(1_A \otimes f) \circ (\alpha \otimes 1_A) = (1_A \otimes g) \circ (\alpha
\otimes 1_A)$ where $f$ and $g$ denote the left and right side of
the equation $(\star)$ respectively. Composing with 
$\varepsilon \otimes 1_A \otimes 1_H$, we get $f=g$ and thus $A$
is a $H$-Hopf module.
 If $H$ has a bijective antipode we can apply the fundamental
theorem of Hopf modules (\cite{[LS],[Mo]}): 
$A \cong A^{coH} \otimes H$ where $A^{coH}$
denotes the invariants of $A$ under the coaction of
$H$. We have $A^{coH} \cong k$ since $A^{coA} \cong k$ and hence
$A \cong H$. The assumption i) implies that $\pi$ is surjective and if $A$
is finite-dimensional, it is also bijective. \square

\medskip

We now come to measured algebras:

\begin{defi}
Let $Z$ be an algebra. A {\rm measure} on $Z$ is a linear form
$\phi : Z \longrightarrow k$ such that the following condition holds:

Let $a \in A$, if $\phi(ab) = 0 \ \forall b \in Z$, then $a=0$.

A {\rm measured algebra} is a pair $(Z,\phi)$ where $Z$ is an algebra
and $\phi$ is a measure on $Z$.
The category of measured algebras (denoted $\mathcal M alg$) as all
algebra morphisms as arrows. The category of finite-dimensional
measured algebras is denoted by $\mathcal M alg_f$.
\end{defi}

\begin{rema}
{\rm
1) Following the classical philosophy for quantum groups and 
quantum spaces, a measure on an algebra is thought as an integration
of functions on the underlying quantum space.

\noindent
2) Let $\phi$ be a linear form on a finite-dimensional algebra $Z$
and let $a \in Z$. We denote by $\phi(a-)$ the linear form on $Z$
defined by $\phi(a-)(b) = \phi(ab)$. Then $\phi$ is a measure
if and only if the map $A \longrightarrow A^*$, $a \mapsto \phi(a-)$
is an isomorphism. This isomorphism is a right $A$-module isomorphism
and hence a measured algebra is a Frobenius algebra. Conversely
a Frobenius algebra always has a measure : the categories
$\mathcal M alg_f$ and $\mathcal Fr$ (Frobenius algebras) are equivalent.

\noindent
3) Let $(Z, \phi)$ be a measured algebra. If $a \in Z$ is invertible,
then $\phi(a-)$ is a measure on $Z$.}
\end{rema}

The next lemma 2.6 will be useful. Before we need:

\begin{lemm}
Let $Z$ be a finite-dimensional algebra. We assume that the base field $k$
contains at least $n={\rm dim} Z$ non-zero distinct elements.
Then every  element of $Z$ is a linear combination of invertible elements.
\end{lemm}

\noindent
{\bf Proof}. We first consider the matrix algebra $M_n(k)$.
If $a \in M_n(k)$ is not invertible, there exists $\lambda \in k^*$
such that $a - \lambda id$ is invertible (by the assumption on the field).
Hence $a = (a - \lambda id ) + \lambda id$ is a linear combination of invertible elements.
Now let us consider the regular representation $\L : Z \longrightarrow
End(Z)$ ($End(Z) \cong M_n(k)$) : $L(a)(b) = ab$. If $L(a)$ is not invertible
there exists $\lambda \in k^*$ such that $L(a) - \lambda id = L(a - \lambda 1)$ 
is invertible. But $L(a - \lambda 1)$ is a right $A$-module map, and hence its 
inverse is also a right $A$-module map, ie is of the form $L(b)$
for some $b \in A$. Thus $a = (a - \lambda 1) + \lambda 1$ is a linear 
combination of invertible elements. \square

\medskip

C. Cibils showed me examples of algebras where the conclusion of the 
lemma is false if the base field does not have enough elements.
In our setting the result over an arbitrary field would be sufficient
for Hopf algebras. 

\begin{lemm}
Let $(Z,\phi)$ be a finite-dimensional measured algebra.
Let us assume that the base field satisfies the assumption of
lemma 2.5. Every linear form
on $Z$ is a linear combination of measures.
\end{lemm}

\noindent
{\bf Proof}.
Let $f \in Z^*$. Then by the assumption $f = \phi(a-)$
for some $a\in Z$. By lemma 2.5 $a$ can be written as a linear combination
of invertible elements, and the result follows from remark 2.5.3. \square

\begin{defi}
Let $A$ be a Hopf algebra. A {\rm measured (right) $A$-comodule algebra}
is a pair $(Z,\phi)$ where $Z$ is a right $A$-comodule algebra
and $\phi$ is a right $A$-colinear measure on $Z$.
The category of measured $A$-comodule algebras (resp.
finite-dimensional measured $A$-comodule algebras) will
be denoted by $\mathcal Mcoalg(A)$ (resp. $\mathcal M coalg_f(A)$):
the morphisms are $A$-colinear algebra maps.
\end{defi}

\begin{ex}
{\rm 
This example is of vital importance to us.
Let $A$ be a finite-dimensional Hopf algebra. Then by \cite{[LS]}, theorem 1.3,
there is a measure $J : A \longrightarrow k$ (the Haar measure)
such that $(A,J)$ is a measured $A$-comodule algebra. Furthermore
$A$ is cosemisimple if and only if $J(1) \not = 0$.}
\end{ex}

\section{Reconstruction results}

Let $\mathcal C$ be a (small) category and let 
$F : \mathcal C \longrightarrow \mathcal Alg_f$ 
(finite-dimensional algebras) be a functor. We are looking for an
analogue of {\bf (1)} in the introduction. It seems that a version
of the forthcoming theorem 3.1 is only available for bialgebras. So we work
with measured algebras, a very natural framework in view of example 2.8.

\begin{theo}
Let $\mathcal C$ be a (small) category and let 
$F : \mathcal C \longrightarrow \mathcal M alg_f$ be a functor.
Then there is a Hopf algebra $A_{aut}(F)$ such that 
$F$ factorizes  through a functor 
$\overline F : \mathcal C \longrightarrow \mathcal M coalg_f(A_{aut}(F))$
followed by the forgetful functor:

$$\xymatrix{
F : \mathcal C \ar[rr] \ar@{-->}[rd]_{\overline F} && \mathcal Malg_f \\
& \mathcal M coalg_f(A_{aut}(F))\ar[ur]}$$

The Hopf algebra $A_{aut}(F)$ has the following universal property:

If $B$ is a Hopf algebra such that $F$ factorizes through
$\mathcal Mcoalg_f(B)$, there is a unique Hopf algebra morphism
$\pi : A_{aut}(F) \longrightarrow B$ such that the following diagram commutes:

$$\xymatrix{
F \ar[rr]^{\alpha} \ar[rd]_{\beta} && 
F \otimes A_{aut}(F) \ar[dl]^{1_F \otimes \pi} \\
& F \otimes B }$$
where $\alpha$ and $\beta$ denote the coactions of $A_{aut}(F)$ and 
$B$ respectively.
\end{theo}

\noindent
{\bf Proof}. We begin by the construction of the Hopf algebra $A_{aut}(F)$.
For this purpose we construct an autonomous monoidal category endowed with a 
fibre functor and we apply Tannaka duality.

Let $\mathcal C(F)$ be the following category. The objects of 
$\mathcal C(F)$ are the tensor products $F(X_1) \otimes ...  \otimes F(X_n)$
with $X_1, ... , X_n \in {\rm ob}(\mathcal C)$ ($k$ is an object of 
$\mathcal C(F)$).
The arrows of $\mathcal C(F)$ are defined to be the linear combinations of
compositions of tensor products of ``elementary arrows'':

- arrows of the form $F(f)$, where $f$ is an arrow in $\mathcal C$;

- $m_X : F(X) \otimes F(X) \longrightarrow F(X)$,
where $X\in {\rm ob}(\mathcal C)$ and $m_X$ is the multiplication of the 
algebra $F(X)$; 

- $u_X : k \longrightarrow F(X)$, where $u_X$ is the unit of the algebra 
$F(X)$;

- $\phi_X : F(X) \longrightarrow k$ where is $\phi_X$ is the measure
associated with the algebra $F(X)$;

- $B'_X : k \rightarrow F(X) \otimes F(X)$, where $B'_X$ is the only
arrow such that $(F(X), B_X, B'_X)$ is a left dual for $F(X)$, see
\cite{[JS]}, sect. 9 ($B_X = \phi_X \circ m_X$ and $B'_X$ exists since
$\phi_X$ is a measure).

It is clear from its definition that $\mathcal C(F)$ is an autonomous
monoidal category and we have a forgetful (fibre) monoidal functor
$U : \mathcal C (F) \longrightarrow Vect_f(k)$ (finite-dimensional 
vector spaces). We are in position to apply Tannaka duality :
let us define $A_{aut}(F) : = {\rm End}^{\vee}(U) = coend(U)$
(see \cite{[JS],[Sc]}). Since $\mathcal C(F)$ is autonomous monoidal,
it follows that $A_{aut}(F)$ is a Hopf algebra.
We know from \cite{[JS],[Sc]} that 
$U : \mathcal C (F) \longrightarrow Vect_f(k)$ factorizes through
$Co_f({\rm End}^{\vee}(U)$). Thus every object of $\mathcal C(F)$
carries a natural $A_{aut}(F)$-comodule structure.
Furthermore for every object $X$ of $\mathcal C$, the maps 
$m_X, u_X,$ and $\phi_X$ are $A_{aut}(F)$-comodule maps, and for 
every arrow $f$ in $\mathcal C$, $F(f)$ is a comodule
map and an algebra map : it follows that we get
the desired functor
$\overline F : \mathcal C \longrightarrow \mathcal M coalg_f(A_{aut}(F))$.

Now let $B$ be a Hopf algebra such that $F$ factorizes through
$\tilde F : \mathcal C \longrightarrow \mathcal M coalg_f(B)$.
It is easy to see that one can extend the coactions of $B$ to all the
objects in $\mathcal C(F)$, and then
$U : \mathcal C(F) \longrightarrow Vect_f(k)$ factorizes through
$\tilde U : \mathcal C(F) \longrightarrow Co_f(B)$ (the maps 
$B'_X$ are automatically comodule maps since the $B_X$ are). 
The universal property of End$^{\vee}(U)$ (cf \cite{[JS],[Sc]}) gives 
the claimed universal property of  $A_{aut}(F)$ . \square

\begin{rema}
{\rm It is clear from the proof that theorem 3.1 also holds when
$\mathcal M alg_f$ is the category of finite-dimensional 
(ie having a dual) measured algebras in a cocomplete braided monoidal
category $\mathcal V$: see \cite{[Sc],[Ly]}.}
\end{rema}

The Hopf algebra $A_{aut}(F)$ is a quantum automorphism group of
the functor $F$ in the sense of \cite{[Wa]}, and as in \cite{[Wa]} we have: 

\begin{prop}
Let $F : \mathcal C \longrightarrow \mathcal M alg_f$ be a functor and
let $A_{aut}(F)$ be the Hopf algebra above.
Then {\rm Hom}$_{k-alg}(A_{aut}(F)) \cong Aut(F)$.
\end{prop}

\noindent
{\bf Proof}. We use the notations of the proof of 3.1.
We have $Aut^{\otimes}(U) \cong {\rm Hom}_{k-alg}({\rm End}^{\vee}(U)$
(\cite{[JS]}) where $Aut^{\otimes}(U)$ is the group of automorphisms
of the monoidal functor $U$. Therefore it is sufficient to prove that
$Aut^{\otimes}(U) \cong Aut(F)$. This quite easy verification is left to the
reader. \square

\begin{ex}
{\rm 1) Let $\mathcal C = \{*\}$ be the category with one object and one
arrow and let $(Z,\phi)$ be a finite-dimensional algebra.
Let $F$ be the functor defined by $F(\{*\}) = (Z,\phi)$.
Let us denote by $A_{aut}(Z,\phi)$ the algebra associated by theorem
3.1: we get algebraic versions of the Hopf $C^*$-algebras constructed in
\cite{[Wa]}. 
See \cite{[Ba]} for the representation theory of $A_{aut}(Z,\phi)$
when $Z$ is a semisimple algebra and $\phi$ is a ``good'' trace on
$Z$.

\noindent
2) in \cite{[Bi]} we constructed the (compact) quantum automorphism
group of  a finite graph. This construction can also be recaptured
algebraically from theorem 3.1.
Let $\mathcal G = (V,E)$ be a finite graph where $V$ is the set of vertices and
$E \subset V \times V$ is the set of edges. Let $C(V)$
and $C(E)$ be the function algebras on $V$ and $E$ respectively.
Let $\phi_V : C(V) \longrightarrow k$ (resp. $\phi_E : C(E)  \longrightarrow k$)
be the classical normalized integration of functions. The quantum automorphism
group of $\mathcal G$ can be described in the following way. Let 
$\mathcal C$ be the category with two objects $\{0,1\}$ and two arrows
$i,j : 0 \rightrightarrows 1$. We define a functor $F : \mathcal C 
\longrightarrow Malg_f$ by
$$F(0) = (C(V),\phi_V) \ ,  \ F(1) = (C(E), \phi_E) \ , \  F(i)= s_* \ ,\ F(j) = t_* \ ,$$
where $s_*$ and $t_*$ are the algebra maps induced by the source and target map
respectively.
Then the algebraic version of $A_{aut}(\mathcal G)$ defined in \cite{[Bi]}
and $A_{aut}(F)$ are isomorphic.}  
\end{ex} 

We now will show that theorem 3.1  enables us to reconstruct a finite-dimensional
Hopf algebra from the category of finite-dimensional measured comodule
algebras and the forgetful functor. This is a quantum generalization
of {\bf (2)}.

\begin{theo}
Let $A$ be a finite-dimensional Hopf algebra.
We assume that the base field $k$ contains at least $n = \dim A$
non-zero distinct elements (cf lemma 2.5).  
Let $\mathcal C = \mathcal M coalg_f(A)$ be the category of
finite-dimensional measured $A$-comodule algebras and let 
$\omega : \mathcal C \longrightarrow \mathcal M alg_f$ be the 
forgetful functor. The Hopf algebras $A$ and $A_{aut}(\omega)$ are
isomorphic.
\end{theo}  

\noindent
{\bf Proof}. We want to apply lemma 2.2. For simplicity  we note
$H = A_{aut}(\omega)$. The universal property of
$H$ yields a Hopf algebra map $\pi : H \longrightarrow A$
and a $H$-comodule structure  $\alpha : A \longrightarrow
A \otimes H$ on $A$ such that 
$(1_A \otimes \pi) \circ \alpha = \Delta_A$.
Let $A_0$ be the trivial $A$-comodule algebra whose underlying algebra
is $A$. We must show that
$A_0$ is also trivial as $H$-comodule (ie $\omega(A_0) \cong \omega(A)_0$).
Let $\psi$ be a measure on $A_0$ : then $\psi$
is a morphism in the category $\mathcal C(\omega)$ of the proof of theorem
3.1., and hence an $H$-comodule morphism $A_0 \longrightarrow k$.
By lemma 2.6 every linear form on $A_0$ is a linear combination of measures,
and hence an $H$-comodule map: this means that $A_0$ is trivial as
$H$-comodule. The map $\kappa = (1_A \otimes m) \circ (\Delta_A \otimes 1_A)
: A \otimes A \longrightarrow A_0 \otimes A$ is an $H$-comodule map
since it belongs to $\mathcal C(\omega)$. The antipode of $H$
is clearly bijective since every object of the category $\mathcal C(\omega)$
is self-dual and hence we are in the situation of lemma 2.2:
$\pi$ is an isomorphism. \square  

\begin{rema}
{\rm
The conclusion of theorem 3.5 is true whenever every linear form
on $A$ is a linear combination of measures and characters.
I ignore if this condition is always realized. More generally, under this 
assumption, and since lemma 2.2 only uses arrows, theorem 3.5 should hold
for a ``finite-dimensional'' Hopf algebra in a cocomplete braided monoidal category
with reasonable dimension theory (an epimorphism between object of 
the same dimension is an isomorphism): we have a classification of
Hopf modules in \cite{[Ly]}. } 
\end{rema}

\section{The $C^*$-algebra case}

The base field is now assumed to be the field of complex numbers.
We want analogues of the results of section 3 in a $C^*$-algebra framework.
The ideas are essentially the same as in the previous sections,
so we will be a little concise.

\begin{defi}
Let $Z$ be a $C^*$-algebra. A measure on $Z$ is a positive and
faithful linear form $\phi$ on $Z$ : $\phi(a^* a) > 0$ for 
$a \not = 0$. A measured $C^*$-algebra is a pair
$(Z,\phi)$ where $Z$ is a $C^*$-algebra and $\phi$ is  a measure
on $Z$. The category of measured $C^*$-algebras 
(denoted $\mathcal M C^*$) as all $C^*$-algebras morphisms
(namely $*$-homomorphisms) as arrows.
The category of finite-dimensional measured $C^*$-algebras
is denoted by $\mathcal M C^*_f$.
\end{defi}

It follows from the Cauchy-Schwartz inequality that a measure 
on a $C^*$-algebra is a measure in the sense of definition 2.3.    

\begin{lemm}
Let $Z$ be a finite-dimensional $C^*$-algebra: every linear
form on $Z$ can be written as a linear combination of measures.
\end{lemm}

\noindent
{\bf Proof}. Let $\phi$ be a positive and faithful trace on $Z$.
Every linear form on $Z$ can be written as $\phi(a-)$ for some
$a \in Z$, and the linear form $\phi(a-)$ is positive and faithful
if and only if $a$ is a positive and invertible element of $Z$.
Any element can be written as a linear combination of 
positive and invertible elements, and thus we get the claimed result.
\square

\medskip

We now want to start with a functor
$F : \mathcal C \longrightarrow \mathcal M C^*_f$ and have an analogue
of theorem 3.1: we will get compact quantum groups 
using Woronowicz' Tannaka-Krein duality \cite{[W2]}.

\smallskip

Let us recall that a {\sl compact quantum group} (\cite{[W1],[W3]},
or Woronowicz algebra or Hopf $C^*$-algebra with unit) is a
pair $(A, \Delta)$ where $A$ is a $C^*$-algebra (with unit)
and $\Delta : A \longrightarrow A \otimes A$ is a coassociative
$*$-homomorphism such that the sets
$\Delta(A)(A\otimes 1)$ and $\Delta(A)(1 \otimes A)$ are both
dense in $A\otimes A$. By abuse of notation a compact quantum group
is often identified with its underlying $C^*$-algebra.
A {\sl morphism} between compact quantum groups $A$ and $B$
is a $*$-homomorphism $\pi : A \longrightarrow B$
such that $\Delta \circ \pi = (\pi \otimes \pi) \circ \Delta$.

Given a compact quantum group $A$ there is a canonically defined
Hopf $*$-algebra $A^o$, which is dense in $A$ (\cite{[W3]}).
A representation of $A$ is a comodule of the Hopf algebra $A^o$.

An {\sl action} of a compact quantum group $A$ on a $C^*$-algebra $Z$ is 
a unital $*$-homomorphism $\alpha : Z \longrightarrow Z \otimes A$
such that there is a dense sub-$*$-algebra $Z^°$ of $Z$
for which $\alpha$ restricts to a coaction 
$\alpha : Z^o \longrightarrow Z^o \otimes A^o$ : ie $Z^o$ is a right
$A^o$-comodule algebra.
A $C^*$-algebra endowed with an action of $A$ is called an
$A$-comodule $C^*$-algebra.

\begin{defi}
Let $A$ be a compact quantum group. A measured $A$-comodule $C^*$-algebra
is a measured $C^*$-algebra $(Z,\phi)$ such that $Z$
is an $A$-comodule $C^*$-algebra and $\phi : Z^o \longrightarrow \mathbb C$
is an $A^o$-comodule map. 
The category of finite-dimensional measured $A$-comodule $C^*$-algebras
will be denoted by $\mathcal M C^*_f(A)$ : the morphisms are
$A$-colinear $*$-homomorphisms.
\end{defi}  

\begin{theo}
i) Let $\mathcal C$ be a (small) category and let 
$F : \mathcal C \longrightarrow \mathcal M C^*_f$ be a functor.
Then there is a compact quantum group $A_{aut}(F)$ such that 
$F$ factorizes  through a functor 
$\overline F : \mathcal C \longrightarrow \mathcal M C^*_f(A_{aut}(F))$
followed by the forgetful functor:

$$\xymatrix{
F : \mathcal C \ar[rr] \ar@{-->}[rd]_{\overline F} && \mathcal M C^*_f \\
& \mathcal M C^*_f(A_{aut}(F))\ar[ur]}$$

The compact quantum group $A_{aut}(F)$ has the following universal property:

If $B$ is a compact quantum group such that $F$ factorizes through
$\mathcal M C^*_f(B)$, there is a unique compact quantum group morphism
$\pi : A_{aut}(F) \longrightarrow B$ such that the following diagram commutes:

$$\xymatrix{
F \ar[rr]^{\alpha} \ar[rd]_{\beta} && 
F \otimes A_{aut}(F) \ar[dl]^{1_F \otimes \pi} \\
& F \otimes B }$$
where $\alpha$ and $\beta$ denote the actions of $A_{aut}(F)$ and 
$B$ respectively.

ii) Let $A$ be a finite quantum group 
(a finite-dimensional Hopf $C^*$-algebra),
Let $\mathcal C = \mathcal M C^*_f(A)$ and let $\omega$ be the forgetful
functor. Then the compact quantum groups $A$ and 
$A_{aut}(\omega)$ are isomorphic.
\end{theo}

\noindent
{\bf Proof}. The proof follows the one of theorem 3.1. Let $X$ be an object
of $\mathcal C$ and let $\phi_X$ be the associated measure on $F(X)$.
We define a scalar product on $F(X)$ by
$\langle x, y \rangle = \phi_X(y^*x)$.
We define a category of Hilbert spaces $\mathcal C'(F)$ in the following way.
The objects of $\mathcal C'(F)$ are the Hilbert spaces 
$F(X_1) \otimes ...  \otimes F(X_n)$
with $X_1, ... , X_n \in {\rm ob}(\mathcal C)$.
The ``elementary arrows'' of $\mathcal C'(F)$ are the ones of 
$\mathcal C(F)$ (defined in the proof of theorem 3.1) plus their adjoints,
and the arrows of $\mathcal C'(F)$ are the linear combinations of
compositions and tensor product of elementary ones.
In this way $\mathcal C'(F)$ is a concrete monoidal $W^*$-category
with conjugates, and we can apply theorem 1.3 of \cite{[W2]}.
We get a compact quantum group $A_{aut}(F)$
as a ``universal admissible pair''. Arguing as in the proof of theorem
3.1 and using theorem 1.3 of \cite{[W2]}, we get assertion i).
The assertion ii) is proved in the same way as theorem 3.5:
we use i) and lemma 4.2 (which plays the role of lemma 2.6).
\square

\bigskip

D\'epartement des Sciences Math\'ematiques,
case 051 

Universit\'e Montpellier II

Place Eug\`ene Bataillon, 34095 Montpellier Cedex 5.

E-mail : bichon\char64 math.univ-montp2.fr

\end{document}